\documentclass[2pt,a4paper]{amsart}

\usepackage{latexsym,amssymb,amsthm,graphics,amsmath,amscd}

\usepackage{amssymb}
\usepackage{eucal}
\usepackage{amsmath}
\usepackage{graphicx} 

\theoremstyle{plain}
\newtheorem{theorem}{Theorem}[section]
\newtheorem*{Theorem B}{Theorem B}
\newtheorem*{Theorem A}{Theorem A}

\newtheorem{definition}{Definition}[section]
\newtheorem{example}{Example}[section]
\numberwithin{equation}{section}

\theoremstyle{remark}
\newtheorem{remark}{Remark}[section]
 \numberwithin{equation}{section}

\def\<{\left < }
\def\>{\right >}
\def\({\left ( }
\def\){\right )}

\def\x{{\bf x}}

\def\2{\#_2 }
\def\x{\bf x}

\begin{document}
\setcounter{page}{1}

\title[Designs in compact symmetric spaces]{Designs in compact symmetric spaces and applications of great antipodal sets}

\author[B.-Y. Chen]{Bang-Yen Chen}
\address{Department of Mathematics, Michigan State University, 
East Lansing, MI 48824, USA}
 \email{chenb@msu.edu}

\begin{abstract}  The theory of designs is an important branch of combinatorial mathematics.  It is well-known in the theory of designs that a finite subset of a sphere is a tight spherical 1-design if and only if it is a pair of antipodal points. On the other hand, antipodal sets and 2-number for a Riemannian manifold are introduced by B.-Y. Chen and T. Nagano \cite{CN82} in 1982.
  An antipodal set is called a great antipodal set if its cardinality is equal to the 2-number.  The main purpose of this paper is to provide a survey on important results in  compact symmetric spaces with great antipodal sets as the designs. In the last two sections of this paper, we present some important applications of 2-number and great antipodal sets to topology and group theory.
\vspace{0.2 cm} 
   \end{abstract}

\maketitle

\section{Introduction}\label{S1}

The theory of designs is an important branch of combinatorial mathematics that deals with the existence, construction
and properties of systems of finite sets whose arrangements satisfy generalized concepts of balance
and/or symmetry.

The concept of $t$-designs of points on an $n$-sphere
 $S^{n} \subset \mathbb R^{n+1}$ was initiated  by  P. Delsarte, J. M. Goethals and J. J. Seidel as an analogy with combinatorial $r$-designs \cite{Dels77}. These spherical $t$-designs have
interesting  properties which are related to several important areas of mathematics \cite{Goe79}.
These ideas were later extended  to the notion of $r$-designs in Delsarte spaces by A. Neumaier  in \cite{Neu81}.
The notion of spherical $t$-designs was also extended as spherical designs of harmonic
index $t$. 
E. Bannai, T. Okuda and M. Tagami \cite{BOT15} studied lower bounds and ``tight'' examples of spherical designs of harmonic index $t$.  A theory of designs in compact symmetric spaces of rank one was given by S. G. Hoggar in \cite{Hog82}. 
Moreover, the classification  of tight designs were developed
by E. Bannai and S. G. Hoggar \cite{BHog85,BHog89},  Y. I. Lyubich \cite{Lyu09}   and S. G. Hoggar \cite{Hog82,Hog1984}, among others.

On the other hand, it is not an easy task to study similar theories in compact symmetric spaces of higher rank.  Nevertheless, there are several studies
in theory of designs in Grassmannian manifolds (see, e.g., \cite{BCN02,BBC04,KO20,Roy10}).
Furthermore, some studies of design theories on unitary groups and complex spheres can be found in
\cite{Nak17,RS09,RS14}, among others.

The concept of the antipodal points in spheres was generalized and developed to antipodal sets in compact symmetric spaces by  B.-Y. Chen and T. Nagano in \cite{CN82,CN88}.  Note that every antipodal set for any compact symmetric space is a finite set and the maximum cardinality is called the 2-number. An antipodal set on a compact symmetric space $M$ is called a great antipodal set whenever  its cardinality is equal to the 2-number of $M$.

The main purpose of this paper is to provide a survey on important results in  compact symmetric spaces with great antipodal sets as designs. Furthermore, in the last two sections of this paper, we present some applications of 2-number and great antipodal sets to topology and group theory. This paper can be regarded as a natural  continuation of my earlier articles \cite{C18,C22,C23}.

\section{Spherical designs}\label{S2}

 The spherical codes and  designs on  an $n$-sphere $S^{n}$ are finite subsets. In 1973,  P. Delsarte \cite{Dels73}   unified the theories of codes and designs on association schemes, and gave the upper bounds for codes and the lower bound for designs by applying linear programming for polynomials associated with metric or cometric association schemes.
From code theoretical viewpoint, a spherical code is to find a finite set $X\subset S^{n}$ such that points
on  $X$ are scattered on $S^{n}$ as far as possible. On the other hand, from design theoretical viewpoint, it is to try to find $X$ which globally approximates the sphere $S^{n}$ very well. 

Of course,  ``What does approximate the sphere  $S^{n}$ well mean?'' is an interesting question. For this question, there exists one very reasonable answer given  in 1977 by P. Delsarte, J. M. Goethals and J. J. Seidel defined as follows.

\begin{definition} \cite{Dels77}  For a natural number $t$, a finite subset $X\subset S^{n}$ is called a {\it spherical $t$-design} if we have
$$\frac{1}{|S^{n}|}\int_{\x\in S^{n}} f({\x})d\sigma({\x})=\frac{1}{|X|}\sum_{u\in X}f(u)$$
 for any polynomial  $f({\x})=f(x_{1},\ldots,x_{n+1})$ of degree at most $t$, with the usual integral on the unit sphere, where $\sigma$ is  the normalized Haar measure and $|X|$ denotes the cardinality of $X$.
\end{definition}

This notion of spherical $t$-designs is an analogue of the classical concept of combinatorial $t$-designs studied in the traditional theory of designs in combinatorics which has caught the hearts of many combinatorialists. This research subject has been investigated on its own but also in connection with several other theories such as lattices, association schemes and number theory as well as cubature formulas. Further, sphere packing problems are found to be related with quantum error correcting codes in the mid of 1990s (see \cite{BCN89,Con96,SS98}).
Moreover, there is a very nice analogy between the theories of codes and designs in the frame work of association schemes formulated in 1973 by P. Delsarte \cite{Dels73} and the theory of spherical codes and designs formulated in 1977 by Delsarte, J. M. Goethals and J. J. Seidel
\cite{Dels77}, which are known today as {\it Delsarte's theory} on association schemes and Delsarte theory on spheres, or slightly more broadly, algebraic combinatorics on association schemes or algebraic combinatorics on spheres (see \cite{BB2009} for details). The essential tool in their works is the addition formula for polynomials; polynomials associated with metric or cometric association schemes, or the Gegenbauer polynomials associated with spheres.
In general, the theory of designs can be given on the Delsarte spaces.

It is well-known that the images of  spherical $t$-designs under an orthogonal transformation are also  spherical $t$-designs. Furthermore, it is also well-known that for any given fixed $n$ and $t$, there always exists a spherical $t$-design $X \subset S^{n}$ whenever the cardinality $|X|$ of $X$ is sufficiently large.

The spherical $t$-designs have been extended in 1980 by A. Neumaier \cite{Neu81} to the concept of $r$-designs in Delsarte spaces, which was so named by an analogy with the $Q$-polynomial association schemes of Delsarte (see, e.g.,  \cite{Dels73,Neu81}). 
The concept of spherical $t$-designs is also generalized as spherical designs of harmonic
index $t$. A finite subset $Y$ on $n$-sphere $S^{n}$ is called a {\it spherical design of harmonic index
$t$} if  it satisfied $\sum_{{\x} \in Y}f({\x})=0$ for any $f\in {\rm Harm}_{t}(S^{n})$ (see \cite{BB2009} for more detail). 
For instance, E. Bannai, T. Okuda  and M. Tagami  investigated in \cite{BOT15} lower bounds and studied examples of  tight spherical designs of harmonic index $t$. For  the study of bounds on antipodal spherical designs with few angles, we refer to \cite{XXY21}.
 
A {\it tight design} is a design whose cardinality is equal to the known natural lower bound.
Any pair of antipodal points in a sphere is known to be a tight 1-design.
Recently,  the classification of tight spherical $t$-designs have been a very attractive research subject.

\section{Design theory in compact symmetric spaces}\label{S3}

In this section, we discuss some generalizations of spherical designs and codes to compact symmetric spaces such as real and complex projective spaces and Grassmannian manifolds
(cf. \cite{BB2009,H} for more details).

\vskip.1in
 \noindent (1) {\it Real and complex projective spaces}.
There are several attempts to generalize the theory of spherical codes and designs to other
spaces. Before Delsarte, Goethals and Seidel studied spherical designs and codes in \cite{Dels77}, they also considered
systems of lines  in $\mathbb R^{n}$ and $\mathbb C^{n}$ (see  \cite{Dels75}). This study is essentially equivalent to considering points in real and
complex projective spaces in which Jacobi polynomials appeared  instead of the Gegenbauer polynomials.
Their attempts tried to find a common framework to study finite subsets in these projective spaces
 and subsets in certain association schemes. They also introduced
several notions, such as Delsarte spaces (see  \cite{Neu80,Neu81,God93}) and polynomial spaces (see \cite{Le92,Le98,Le98.2}). 
The similarities of the theories between these continuous spaces and association schemes provided evidence from  of Delsarte's 1973 works \cite{Dels73} and the 1977 works of
Delsarte, Goethals and Seidel  given in \cite{Dels77}. 
The most natural set up for the continuous spaces are  compact symmetric spaces of rank one.
This theory was first developed  in that framework  by S. G. Hoggar \cite{Hog82} in 1982. Furthermore, in the theory of codes and designs, the concept of tight $t$-designs and the classification problems of tight $t$-designs  were developed in a series of papers by E.~Bannai and S. G.~Hoggar in \cite{BHog85,BB2009,Hog1984,Hog1989,Ho90}.

\vskip.1in
\noindent (2) {\it  Grassmannian manifolds.}
The Grassmannian manifolds are natural generalization of projective spaces in which one-dimensional subspaces are considered. As a well-known fact, we may
regarded Grassmannian manifolds $G^{\mathbb F}(n,k)$ as the set of $k$-dimensional subspaces of an $n$-dimensional
vector space $V$ over the field $\mathbb F=\mathbb R, \mathbb C$ or $\mathbb H$. Obviously, Grassmannian manifolds are compact symmetric spaces of higher ranks if $k>1$.

The concept of $t$-designs in Grassmannian manifolds was first introduced in 2002 by C. Bachoc, R. Coulangeon and G. Nebe in  \cite{BCN02}. Basic properties $t$-designs in Grassmannian manifolds were also studied in \cite{BCN02}. 
Some types of Delsarte theory as well as the formulation of tight $t$-designs
was further investigated in 2004 by C.~Bachoc, E.~Bannai and R.~Coulangeon  \cite{BBC04}. On the other hand, the classification problems of tight
$t$-designs remain open til now. 
For further works on algebraic combinatorics of Grassmannian manifolds, we refer to \cite{Bac05,Bac06,BN06}, among others.

\section{Spherical 1-designs, great antipodal sets and 2-number}\label{S4}

The following results are well-known in spherical design theory.
Let $Y$ be a finite subset of an $n$-sphere $S^{n}$. Then we have:

\begin{itemize}
\item[(1)] If $Y$ is a spherical 1-design (or equivalently, a spherical design of harmonic index 1),
then the cardinality of $Y$ satisfies $|Y| \geq 2$.

\item[(2)] $Y$ is a spherical 1-design with $|Y| = 2$ if and only if $Y$ is a pair of antipodal points.

\item[(3)]  Let  $T$ be a set of integers $t$ such that a pair of antipodal points on $S^{n}$ is a spherical
design of harmonic index $t$. Then $T$ consists of all odd natural numbers.
\end{itemize}

Recall the following definition from \cite{CN82,CN88}.

\begin{definition} {\rm An antipodal set $A(M)$ of a Riemannian manifold $M$ is a discrete subset of $M$
such that any two points in $A(M)$ are antipodal on some closed geodesic of $M$.
A {\it maximal antipodal set} of $M$ is an antipodal set which is not a proper subset
of any other antipodal set.}
\end{definition}

Note that the cardinalities of antipodal sets are bounded for each compact Riemannian manifold (for a clear proof, see  \cite{TT13}).
The notions of maximal antipodal sets allow us to define ``2-number'', denoted by $\#_{2}M$, which
is the largest cardinal number among all maximal antipodal sets. A maximal
antipodal set $A(M)$ is called a {\it great antipodal set} if $A(M)$ satisfies $\#(A(M)) = \#_{2}M,$ where $\#(A(M))$ denotes the cardinal number of $A(M)$.

For 2-number of Riemannian manifolds, we have \cite{CN88}

(a) $\#_{2}M\geq 2$ for every compact Riemannian manifold $M$;

(b) $\#_{2} S^{n}= 2$ and $\#_{2}RP^{n} = n + 1$.

\begin{remark} Statement (a) follows from the fact that every compact Riemannian manifold admits a closed geodesic (see \cite{LF51}).
\end{remark}

The notion of symmetric $R$-spaces were defined in 1965 independently by Nagano \cite{Na65} and Takeuchi \cite{Ta65}  as $R$-spaces  which are at the same time compact symmetric spaces. Hence, a  symmetric $R$-spaces admits a transitive action of a center-free non-compact semi-simple Lie group and the corresponding stabilizer of a point is a certain maximal parabolic subgroup. But the name ``{\it symmetric R-space}'' was coined by M. Takeuchi in \cite{Ta65}.

The family of symmetric $R$-spaces have been completely classified by S. Kobayashi and T. Nagano  in \cite{KN64} which consists of (up to minor modifications)
\begin{itemize}
  \item[{\rm (a)}] All Hermitian symmetric spaces of compact type.

 \item[{\rm (b)}]  Grassmann manifolds $O(p+q)/O(p)\times O(q), Sp(p+q)/Sp(p)\times Sp(q)$.

 \item[{\rm (c)}]  The classical groups $SO(m),\,U(m),\,Sp(m)$.

 \item[{\rm (d)}]  $U(2m)/Sp(m),\, U(m)/O(m)$.

 \item[{\rm (e)}]  $(SO(p+1)\times  SO(q+1))/S(O(p)\times O(q))$, where $S(O(p)\times O(q))$ is  the subgroup  of $SO(p+1)\times SO(q+1)$ consisting of matrices of the form:
$$\begin{pmatrix} \epsilon& 0 & & \\ 0 & A&&\\ &&\epsilon& 0\\
&&0&B \end{pmatrix},\quad \epsilon=\pm1,\quad A\in
O(p),\quad B\in O(q),$$
 (This symmetric $R$-space is covered twice by $S^p\times S^q$).

 \item[{\rm (f)}]  The Cayley projective plane $FII={\mathcal O}P^2$. 

 \item[{\rm (g)}]  The three exceptional spaces $E_6/Spin(10)\times S^{1}, E_7/E_6\times S^{1},$ and $E_6/F_4.$
\end{itemize}

\vskip.1in
Great antipodal sets and the 2-number on symmetric $R$-spaces have several nice properties; for examples, we have:

\begin{itemize}
  \item[{\rm (i)}] Great antipodal sets on a symmetric $R$-space are unique up to congruence
(see \cite{TaT12,TT13}).

  \item[{\rm (ii)}]  Great antipodal sets on a symmetric $R$-space are certain orbits
of Weyl groups  (see \cite{Ta89}).
\end{itemize}

\section{Unitary groups with great antipodal sets as designs}\label{S5}

The compact simple Lie group $U(n)$ is a compact symmetric space (with respect to a bi-invariant metric) and it has the point-symmetry at each given
point in $U(n)$. A great antipodal set on $U(n)$ is
an analogue of a pair of antipodal points on spheres which is a tight spherical 1-design.  

Let $G$ be a compact group together with a closed subgroup $K$. Put $M = G/K$ and let $\mu$  be a positive finite $G$-invariant measure on $M$. Suppose $\rho$ is the natural unitary representation of $G$ on the Hilbert space $L^{2}(M,\mu)$. Choose a $\rho(G)$-invariant irreducible closed subspace
$\mathcal H$ of $L^{2}(M,\mu)$.
Denote by  ${\mathcal C}(U(n))$  the space of $\mathbb C$-valued continuous functions on $U(n)$. Then  
${\mathcal C}(U(n))$ is  a $G$-space as $$G \times {\mathcal C}(U(n)) \to {\mathcal C}(U(n)); \;\;(x_{1},x_{2}) \cdot f (x) = f(x_{1}^{-1} x x_{2}).$$
Let ${\mathcal H}_{\lambda}$ be a $G$-subspace of ${\mathcal C}(U(n))$  which is isomorphic to 
$\rho_{\lambda}\boxtimes \rho_{\lambda}^{*}$. Then it is known
that the structure of ${\mathcal H}_{\lambda}$ in ${\mathcal C}(U(n))$  exists uniquely, i.e.,  if $W \subset {\mathcal C}(U(n))$ is a
subrepresentation which is isomorphic to $\rho_{\lambda}\boxtimes \rho_{\lambda}^{*}$, then we have $W = {\mathcal H}_{\lambda}$ (see \cite{Ku21} for more details).

In 2001, H. Kurihara  \cite{Ku21} investigated the relationship between great antipodal sets on $U(n)$ and design theory on $U(n)$.  He established a nice link between great antipodal sets on $U(n)$ and designs. He also obtained a nice link between great antipodal sets on $U(n)$ and the Hamming cube graph ${\mathcal Q}_{n}$.

\begin{definition} {\rm Fix a subset ${\mathcal T}\subset (\mathbb Z^{n})_{+}$ and let $Y$ be a finite subset of $U(n)$. Then $Y$ is called a {\it ${\mathcal T}$-design} if the following formula:
$$\int_{U(n)} f({\x}) d\sigma({\x})=\frac{1}{|Y|}\sum_{x\in Y} f(\x)$$
is satisfied for every $f\in \oplus_{\lambda\in {\mathcal T}} {\mathcal H}_{\lambda}$.  When the subset ${\mathcal T}$ is a singleton $\lambda$, then $Y$ is called a {\it $\lambda$-design}.}
\end{definition}

For $\lambda = (\lambda_{1},\ldots,\lambda_{n}) \in (\mathbb Z^{n})_{+}$, let us put ${\rm diff}(\lambda) := \lambda_{1}-\lambda_{n}$. Note that ${\rm diff}(\lambda)$  is invariant
for the equivalence relation $\sim$. Thus we may define ${\rm diff}([\lambda])$.

H. Kurihara {\rm \cite{Ku21}} proved the following result for a great antipodal set $S$ on $U(n)$.

\begin{theorem} \label{T:5.1} For a great antipodal set $S$ on $U(n)$,
\begin{itemize}
  \item[{\rm (1)}] 
   If $[\lambda]$ is odd, then S is a $[\lambda]$-design.
 
   \item[{\rm (2)}]  If $[\lambda]$ is even and S is a $[\lambda]$-design, then ${\rm diff}([\lambda]) < n - 1$. In particular, there
exist only finitely many even $[\lambda]$ such that $S$ is a $[\lambda]$-design.
   \end{itemize}
\end{theorem}

\begin{remark} Theorem \ref{T:5.1} shows that  great antipodal sets are  ``good'' designs  on $U(n)$. 
\end{remark}

H. Kurihara also provided the following examples in {\rm \cite{Ku21}}.

\begin{example} {\rm For small $n$,  the condition on $[\lambda]$ carries that a great antipodal set
$S$ is a $[\lambda]$-design.}
 
\begin{itemize}
  \item[{\rm (1)}]   {\rm When $n = 2$,  $S$ on $U(2)$ is an even $[\lambda]$-design if and only
if $[\lambda]= [(1, 1)]$.}
 
 \item[{\rm (2)}]   {\rm When $n = 3$, $S$ on $U(3)$ is an even $[\lambda]$-design if and only if
$[\lambda]= [(1, 1, 0)]$ or $[\lambda]=[(2, 1, 1)]$.}
 
 \item[{\rm (3)}]   {\rm $S$ on $U(4)$ is an even $[\lambda]$-design if and only if the following six classes: 
$[\lambda]=$ \\
$[(1, 1, 0, 0)], [(2, 1, 1, 0)], [(1, 1, 1, 1)], [(3, 1, 1, 1)], [(2, 2, 1, 1)], [(3, 3, 3, 1)]$.}

 \item[{\rm (4)}]  {\rm $S$ on $U(5)$ is an even $[\lambda]$-design if and only if the following 12 classes: 
$[\lambda]=$ \\
$[(1, 1, 0, 0, 0)], [(2, 1, 1, 0, 0)],  [(1, 1, 1, 1, 0)],  [(3, 1, 1, 1, 0)],  [(2, 2, 1, 1, 0)]$,\\
 \;  $[(3, 3, 3, 3, 0)],  [(2, 1, 1, 1, 1)], \, [(4, 1, 1, 1, 1)],  (3, 2, 1, 1, 1)],  [(2, 2, 2, 1, 1)],$\\
 $[(3, 3, 3, 2, 1)], \; [(4, 3, 3, 3, 1)]$.}

\end{itemize}
\end{example}

Let $[n]$ denote the set $\{1, 2, . . . , n\}$ and $Q$ be the $n$-ary Cartesian product of the two elements set $\{1,-1\}$. The {\it Hamming cube graph} $\mathcal Q_{n}$ of degree $n$ is the graph with the vertex set $Q$ and two vertices are adjacent if they differ in precisely one coordinate.

Let $g$ be a bi-invariant metric on $U(n)$ and let $${\rm dist}_{g} : U(n)\times U(n)\to \mathbb R$$
be the distance function on $U(n)$ via $g.$ For a great antipodal set $S$ of $U(n)$, we
define the {\it minimum distance} of $S$ with respect to $g$ by
$${\rm md}_{g}(S) := \min \{{\rm dist}_{g}(x, y) : x,y \in S, x\ne y\}.$$
Now, consider the graph defined as follows: The vertices set is $S$ and two vertices $x$ and $y$ in $S$
are adjacent if we have ${\rm dist}_{g}(x, y) = {\rm md}_{g}(S)$. Let $E_{g}$ denotes the edge set of this graph.

In {\rm \cite{Ku21}}, H. Kurihara  also proved the following  link between great antipodal sets on $U(n)$ and the Hamming cube graph ${\mathcal Q}_{n}$.

\begin{theorem} \label{T:5.2} For any bi-invariant metric $g$ on $U(n)$, the graph $(S, E_{g})$ is a Hamming cube $\mathcal Q_{n}$.
\end{theorem}

\begin{remark} For further results in this direction, see H. Kurihara \cite{Ku21}. For tight 4-designs in Hamming association schemes, we refer to A. Gavrilyuk, S. Suda and J. Vidali \cite{GSV20}.
\end{remark}

\section{Grassmannian manifolds with great antipodal sets as designs}\label{S6}

In this section we present some results  on characterizations of great antipodal sets as design in complex Grassmannian manifolds among certain designs with the smallest cardinalities. 
Note that great antipodal sets are in differential geometry. On the other hand, the theory of designs is related to algebraic combinatorics or representation theory.

For a 1-design $X$ on a complex projective $n$-space ${\mathbb C}P^{n}$, we have the following inequality: $$|X|\geq n+1,$$ and $X$ is tight if $|X| = n+1$ (see  \cite{Hog82}). 

It is also known that for a complex projective space, great antipodal sets can be characterized by tight 1-designs. Therefore, we have the following result \cite[Fact 3.1, page 440]{KO20}.

\begin{theorem}\label{T:6.1}
Let $S$ be a finite subset of a complex projective space. Then the following two statements are equivalent:
\begin{itemize}
  \item[{\rm (a)}]  $S$ is a great antipodal set.
   \item[{\rm (b)}] $S$ is a tight 1-design.
\end{itemize}
\end{theorem}

In order to state the main results of this section, we put 
\begin{align}
&{\mathcal E} =\{(\overbrace{1,\ldots,1}^{i},\overbrace{0,\ldots,0}^{m-i}: i=0,1,\ldots,m\},
\\& {\mathcal F} =\{((2,\overbrace{1,\ldots,1}^{i-1},\overbrace{0,\ldots,0}^{m-i}): i=0,1,\ldots,m\}.
\end{align}

For great antipodal sets of a complex Grassmannian manifold,
 H. Kurihara and T. Okuda \cite{KO20} proved the following two main results of this section.

\begin{theorem} \label{T:6.2} A great antipodal set  of a complex Grassmannian manifold $G^{\mathbb C}(m,n)$ is an ${\mathcal E}$-design with the smallest cardinality.
\end{theorem}

\begin{theorem} \label{T:6.3}
Let $S$ be a finite subset of  $G^{\mathbb C}(m,n)$. Then the following conditions are equivalent:
\begin{itemize}
  \item[{\rm (i)}]  $S$ is a great antipodal set on $G^{\mathbb C}(m,n)$.
  
   \item[{\rm (ii)}] $S$ is an ${\mathcal E}\cup {\mathcal F}$-design on $G^{\mathbb C}(m,n)$  with the smallest cardinality.
\end{itemize}
\end{theorem}

\begin{remark} Theorems \ref{T:6.2} and \ref{T:6.3} show that great antipodal sets  are  ``good'' designs for complex Grassmannian manifolds  \cite{KO20}.
\end{remark}

\begin{remark} A great antipodal set $S$ is a 1-design on of a complex Grassmannian manifold according to Theorem \ref{T:6.2}. However, it cannot be a 2-design in general according to \cite[Remark A.1, page 463]{KO20}.
\end{remark}

\section{Cubature formulas for great antipodal sets in Delsarte theory}\label{S7}

Let $\Omega$
 be a subset of a Euclidean $n$-space $\mathbb R^{n}$. Consider the following integral:
$$\int_{\Omega} f(\x)\mu(\x)d\x$$
where  $\mu$ is a positive weight function on $\Omega$, where we assume that all polynomials of up to sufficiently large
degrees are always integrable. We also assume that the weight function $\mu(x)$ is normalized such that
$\int_{\Omega} \mu(\x)d\x=1$ holds.
    
 \begin{definition} {\rm \cite{HOS16,Ok17}} $($Cubature formula of degree $t)$. {\rm Let $X=\{ u_{1},\ldots,u_{N}\}$ be a finite subset of $\Omega\subset \mathbb R^{n}$. Then, for any polynomial $f(\x)$ of degree $t$ over real field $\mathbb R$ and $N$ positive real numbers  $\lambda_{1},\ldots, \lambda_{N}$, the equation:
$$ {\int_{\Omega} f({\x})\mu({\x})d{\x}} =\sum_{i=1}^{N} \lambda_{i}f(u_{i})$$
 is called a {\it cubature formula of degree $t$ with $N$ points}, where $\lambda_{1},\ldots, \lambda_{N}$  are independent of the choice of the polynomial $f(\x)$.}
\end{definition}

T. Okuda and H. Kurihara \cite{Ok17} established a formulation of Delsarte theory for finite subsets of
compact symmetric spaces. As an application, they proved that great antipodal subsets of complex Grassmannian manifolds give rise to the following cubature formulas.

\begin{theorem} {\rm \cite{Ok17}} Let $X$ be a great antipodal set of $G^{\mathbb C}(m,n)$ and $\mathcal{S}$ be  a
finite-dimensional functional space on $G^{\mathbb C}(m,n)$. Then  the
following cubature formula:
\begin{equation}\label{7.1} \frac{1}{{\rm vol}(G^{\mathbb C}(m,n))} \int_{G^{\mathbb C}(m,n)} fd_{\mu_{G^{\mathbb C}(m,n)}}=\frac{1}{|X|}\sum_{x\in X}f(x),\;\;   \forall f\in \mathcal{S},\end{equation}
holds, where $\mu_{G^{\mathbb C}(m,n)}$ is a $U(n)$-invariant Haar measure  and ${\rm vol}(G^{\mathbb C}(m,n))$
is the volume of  $G^{\mathbb C}(m,n)$ with respect to the measure $\mu_{G^{\mathbb C}(m,n)}$. Furthermore, any great antipodal subset $X$ has the minimum cardinality as a finite subset of $G^{\mathbb C}(m,n)$
such that the formula \eqref{7.1} holds
\end{theorem}

\section{Links between Euler number and  2-number}\label{S8}

\subsection{2-number and Euler number}\label{S8.1}

The author and T. Nagano proved  in \cite{CN88} the following links between 2-number and Euler number.

\begin{theorem}\label{T:8.1} $\#_{2}M \geq  \chi(M)$ for a compact symmetric space $M$.
\end{theorem}

\begin{theorem}\label{T:8.2}  $\#_{2}M = \chi(M)$ for a compact Hermitian symmetric space
of semi-simple type.
\end{theorem}

\subsection{2-number and covering map}\label{S8.2}

 The following links between 2-number and covering map were also obtained in \cite{CN88}.
     
\begin{theorem} \label{T:8.3} Let $M$ and $M''$ be compact symmetric spaces. If $M$ is a double covering of $M''$,  then we have $\#_{2}M\leq 2\, \#_2(M'').$
 \end{theorem}
 
\begin{remark} Inequality $\#_{2}M\leq 2\, \#_2(M'')$ in Theorem \ref{T:8.3} is optimal, because the equality  holds for group manifold $SO(2m)$, $m>2$.
\end{remark}

 For $k$-fold coverings with odd $k$, we have
 
\begin{theorem}\label{T:8.4} Let $\phi :M\to N$ be a $k$-fold covering  between compact symmetric spaces.  When $k$ is odd, we have $\#_{2}M=\#_2(N).$
 \end{theorem}

\begin{remark} Recently, M. S.~Tanaka and H.~Tasaki established a refinement of Theorem  \ref{T:8.4} in the case of compact Lie groups, see \cite{TTa24} for details.
\end{remark}

\subsection{2-number and homology}\label{S8.3}

 M. Takeuchi  \cite{Ta89}  proved in 1989 the following link between 2-number and $\mathbb Z_2$-homology for symmetric $R$-spaces. 

\begin{theorem}\label{T:8.5}  
$\,\#_{2}M=\dim H(M,\mathbb Z_2)$  for every symmetric $R$-space $M$, where $H(M,\mathbb Z_{2})$ is the $\mathbb Z_2$-homology group of $M$. \end{theorem}

\section{Applications to group theory}\label{S9}

In 1953, A. Borel and J.-P. Serre  \cite{BS53} defined the 2-rank,  $r_2G$, of a compact Lie group $G$ as the  maximal possible rank of the elementary 2-subgroups of $G$. In 
 \cite{BS53} they also proved the following two results:
\vskip.05in

\begin{itemize}
  \item[{\rm (1)}] 
 $rk(G)\leq r_2(G)\leq 2\, rk(G)$ and 
 
 \item[{\rm (2)}] $G$ has 2-torsion if $rk(G)<r_2(G)$ 
\end{itemize}
for any compact connected Lie group $G$, where $rk(G)$ is the ordinary rank of $G$.
\vskip.05in

Borel and Serre \cite{BS53} were able to determine the 2-rank of simply-connected simple Lie groups $SO(n), Sp(n), U(n), G_{2}$ and $F_{4}$. They also proved that $G_{2}$, $F_{4}$ and $E_{8}$ have 2-torsion. 
On the other hand, they pointed out in  \cite[page 139]{BS53} that they were unable to determine 2-rank for $E_{6}$ and $E_{7}$ among others. 
To settle this Borel--Serre's problem, the author and T. Nagano established in  \cite{CN82} the following simple links between 2-number and 2-rank.
 
 \begin{theorem}\label{T:9.1} For a connected compact Lie group $G$,  we have
$ \#_{2}G=2^{r_{2}G}.$
 \end{theorem}
  
\begin{theorem}\label{T:9.2} \cite{CN82} For two connected compact Lie groups  $G_{1}$ and $G_{2}$,  we have $ \#_{2}(G_{1}\times G_{2})=2^{r_{2}G_{1}+r_{2}G_{2}}.$
\end{theorem}

Since we are able to determine 2-numbers for (most) compact symmetric spaces in \cite{CN88} via our $(M_{+},M_{-})$-theory developed earlier in 1970s--1980s,   
 we are able to settle this problem of Borel and Serre's problem on 2-ranks of  Lie groups (see \cite{CN77,CN78,CN82,CN88} and  \cite{C87,C89,C13,C18,C22}). 
 
 Now, we state the 2-ranks of compact Lie groups as follows.

\subsection{Classical groups}\label{S9.1}

 For classical groups we have:

\begin{theorem}\label{T:9.3}  For the unitary group $U(n)$, we have 
  \begin{equation}\notag r_{2}(U(n)/{\mathbb Z}\mu)= 
 \begin{cases} n+1 & \text{if $\mu$ is even and $n=2$ or $4$;}\\ n & \text{otherwise,}\end{cases} \end{equation}
 where ${\mathbb Z}\mu$ is a cyclic normal subgroup  of order $\mu$.
 \end{theorem}

 \begin{theorem}\label{T:9.4} For $SU(n)$, we have
 \begin{equation}\notag  r_{2}(SU(n)/{\mathbb Z}\mu)= 
 \begin{cases} n+1 & \text{for $(n,\mu)=(4,2)$;}
 \\ n  & \text{for $(n,\mu)=(2,2)$ or $(4,4)$;}
 \\ n-1 & \text{for the other cases.} \end{cases} \end{equation}
 \end{theorem}

 \begin{theorem}\label{T:9.5} One has $r_{2}(SO(n))=n-1$, and for $SO(n)^{*}$ we have
 \begin{equation}\notag  r_{2}(SO(n)^{*})= 
 \begin{cases} 4 & \text{for $n=4$;}\\ n-2 & \text{for $n$ even $>4$}.\end{cases} \end{equation}
 \end{theorem}

 \begin{theorem}\label{T:9.6} We have
 \begin{itemize} 
\item[{\rm (a)}] $r_{2}(O(n))= n$;
 
\item[{\rm (b)}] $r_{2}(O(n)^{*})$ is $n= 2$ or $4$, while it is $n-1$ otherwise,  
\end{itemize} 
where  $O(n)^{*}=O(n)/\{\pm 1\}$.\end{theorem}

 \begin{theorem}\label{T:9.7} We have $r_{2}(Sp(n))=n$, and for $Sp(n)^{*}$ we have
  \begin{equation}\notag  r_{2}(Sp(n)^{*})= \begin{cases} n+2 & \text{for $n=2$ or $4$}\\ n+1 & \text{otherwise}.
\end{cases} \end{equation}
\end{theorem}

 \subsection{Spinors, semi-spinors and $Pin(n)$}\label{S9.2}
 
For spinor we have the following.

 \begin{theorem}\label{T:9.8} We have
 \begin{equation}\notag
  r_{2}(Spin(n)) =\begin{cases} r+1 & \text{if  $\,n \equiv -1,0$ or {\rm 1 (mod 8)}}\\ 
 r &\text{otherwise}.\end{cases}
 \end{equation}  
 where $r= [\frac{n}{2}]$ is the rank of $Spin(n)$.
 \end{theorem}

 \begin{theorem}\label{T:9.9}  We have
$ r_{2}(Spin(n+8))= r_{2}(Spin(n))+4$ for $n\geq 0. $
 \end{theorem}
 
 For semi-spinor group $SO(4m)^{\#}=Spin(4m)/\{1, e_{((4m))}\}$, we have

 \begin{theorem}\label{T:9.10}  Let $r$ be the rank $2m$ of $SO(4m)^{\#}$. Then we have:
 \begin{equation}\notag
  r_{2}(SO(4m)^{\#}) =\begin{cases} 3 & \text{if  $m=1$}\\  6 &\text{if $m=2$,}\\ r+1 &\text{if $m$ is even $>2$, }\\
 r &\text{if $m$ is odd $>1$}.\end{cases}
 \end{equation}
  \end{theorem}

$Pin(n)$ was introduced by M. F. Atiyah, R. Bott and A. Shapiro in \cite{ABS64}. 

 \begin{theorem}\label{T9.11} We have $r_{2}(Pin(n))= r_{2}(Spin(n + 1))$, $n\geq 0$. \end{theorem}

 \subsection{Exceptional groups}\label{S9.3}

For exceptional Lie groups we have

\begin{theorem}\label{T:9.12}  $r_{2}G_{2}=3,\,  r_{2} F_{4}=5,\,  r_{2}E_{6}=6,\, r_{2}E_{7}=7,\, r_{2}E_{8}=9$ and $r_{2}E^{*}_{6}=6$.\end{theorem}

\begin{remark} $r_2 G_2=3$ and $r_2 F_4=5$ above are due to  A. Borel and J.-P. Serre \cite{BS53}.
\end{remark}

\end{document}